\documentclass[a4paper,10pt]{article}
\usepackage{amssymb}
\usepackage[cp1250]{inputenc}
\usepackage{amsfonts}
\usepackage{amsthm}
\usepackage[english]{babel}
\usepackage{amsmath}
\usepackage[T1]{fontenc}
\usepackage{latexsym}

\usepackage{color}

\newtheorem{theorem}{Theorem}[section]

\newtheorem{proposition}[theorem]{Proposition}
\newtheorem{corollary}[theorem]{Corollary}
\theoremstyle{definition}
\newtheorem{remark}[theorem]{Remark}

\newcommand{\zfrac}[2]{#1/#2}

\begin{document}

\author{Jacek
Jakubowski and
 Maciej Wi\'sniewolski }\thanks{Research partially supported by Polish MNiSW grant N N201 547838.}

\title{On some Brownian functionals and their applications to
moments in lognormal and Stein stochastic volatility models}

\maketitle

\begin{center}
{\small
 Institute of Mathematics, University
of Warsaw \\
  Banacha 2, 02-097 Warszawa, Poland \\
e-mail: {\tt jakub@mimuw.edu.pl } \\
and \\
  {\tt wisniewolski@mimuw.edu.pl } }
\end{center}

\begin{abstract}
The aim of this paper is to present the new results concerning some
functionals of Brownian motion with drift and present their
applications in financial mathematics. We find a  probabilistic
representation of the Laplace transform of special functional of
geometric Brownian motion  using the squared Bessel and radial
Ornstein-Uhlenbeck processes. Knowing the transition density
functions of the above we obtain computable formulas for certain
expectations of the concerned functional. As an example we find the
moments of processes  representing an asset price in the lognormal
volatility ans Stein models.
 We also present links among the  geometric Brownian motion,
the Markov processes studied by Matsumoto and Yor and the hyperbolic
Bessel processes.
\end{abstract}

\noindent
\begin{quote}
\noindent  \textbf{Key words}: geometric Brownian motion,
Ornstein-Uhlenbeck process, Laplace'a transform,  Bessel process,
hyperbolic Bessel process

\ \\
\textbf{2010 AMS Subject Classification}: 91B25, 91G20, 91G80,
60H30.

\textbf{JEL Classification Numbers}: G12, G13.
\end{quote}

\section{Introduction}
The aim of this paper is to present the new results concerning some
 functionals of Brownian motion with drift and their
applications to financial mathematics.  The laws of many different
functionals of Brownian motion have been studied in recent years
(see, among others, \cite{DYM}, \cite{SB}, \cite{DGY}, \cite{MatII},
\cite{MatIII}), but some of the obtained results can not be
effective used in application.
 The distribution of
$\int_0^te^{B_u^{(\mu)}}du$, where $B_t^{(\mu)}=B_t + \mu t$ with a
standard Brownian motion  $B$, is an example of such situation. This
distribution can be characterized by Hartman-Watson distribution,
but the oscillating nature of the last causes the difficulties in
numerical calculations connected with this functional (see
\cite{BRY}).
  We study the laws of special
functionals of geometric Brownian motion, and find results
convenient for numerical applications.
 We investigate the functionals of geometric Brownian motion $Y^{(\mu)}_t := \exp \Big(B_t +
\mu t\Big) $ for $\mu \in R$. In particular we study properties of the functionals
$ \Gamma_t = \frac{Y_t^{(-1/2)}}{1+\beta\int_0^tY_s^{(-1/2)}ds}$ and
$1+ \beta A^{(\mu)}_t$ for  $\beta>0$, where $A^{(\mu)}_t
:= \int_0^t(Y^{(\mu)}_u)^2du$.
 We deliver the probabilistic representation for Laplace transform of
$\Gamma$.
 In our probabilistic representation of
Laplace transform of $\Gamma$ we use the  squared Bessel and radial
Ornstein-Uhlenbeck processes. Knowing the transition density
functions of these processes we obtain computable formulas for
certain expectations of the concerned functionals. One of the many
advantages of the new result is the fact that they can be
effectively used in numerical computations. As an example we compute
the moments $\mathbb{E}X_t^{\alpha}$, for $\alpha > 0$, of the
processes $X_t$ representing an asset price in an important
stochastic volatility model - in the lognormal volatility model. The
necessity of computing moments results from the problems of pricing
derivatives (for instance, the broad class of interest rate
derivatives necessity  a "convexity correction" to the forward rate
price; for details see e.g. \cite{BM}) as well as from the need of
approximations of characteristic functions of random variables with
very complicated distributions.

We now give a detailed plan of  this paper. In subsection 2.1 we
present a method of calculating the moments $
\mathbb{E}\Gamma_t^{k}$ for $k\in \mathbb{Z}$ (Proposition
\ref{recurence}, Remarks \ref{uw1}, Corollary \ref{wn2}) and
investigate the connection of functional $\Gamma$ with a  hyperbolic
Bessel process (Theorem \ref{tw:rephb}). The general connections
between hyperbolic Bessel processes and functionals of geometric
Brownian motion are presented in  \cite{jw}.
 In subsection 2.2 we investigate the  different properties of
$1+ \beta A^{(\mu)}_t$. We find two different probabilistic
representations of the Laplace transform of $(1+\beta
A^{(\mu)}_t)^{-1}$ (Theorems \ref{tw:qu} and \ref{tw:freeb}), the
form of $ \mathbb{E}\ln (1+\beta A^{(\mu)}_t)$ and $ \mathbb{E}
(1+\beta A^{(\mu)}_t)^{-1}$ (Theorem \ref{1st}).
  It turns out
 that for an arbitrary strictly positive random variable we can find a
representation of the Laplace'a transform of $(s+\xi)^{-1}$ for
$s\geq 0$, in terms of a  squared Bessel process (Theorem
\ref{tw:og-rep}). Moreover, we find
 some interesting connections between $\mathbb{E}((1 +\beta
A_t^{(\mu)})^{-1})$ and the conditional expectation of functionals
of geometric Brownian motion with opposite drift.
 Notice that we establish all results for a fixed
$t$.  Section 3 is an illustration of using the previous results
in mathematical finance. We assume that the asset price process $X$
satisfies
 $dX_t = Y_tX_tdW_t$ with $Y$  being a (GBM) (this model is
called the lognormal stochastic volatility model or the Hull-White
model, see \cite{Hul}) and $Y$ being an Ornstein-Uhlenbeck process
(OU) (the Stein model, see \cite{SS}).
 The distribution of  the asset
price for the lognormal stochastic volatility model is known  but
degree of complication and numerical obstacles encourage  to look
for simpler approximations.
 Jourdain
\cite{Jou} has given conditions on existence of the moments, if $Y$
is a (GBM), but not mentioned about how to compute it.
 In this work we
find that the moment is equal to the Laplace'a transform of  the
process $\Gamma$ (Theorem \ref{tw:momrep1}).  We can also express moments of order
$\alpha>1$ in terms of the hyperbolic Bessel process (Theorem
\ref{tw:repr-Bess}). For the model with random time $T_\lambda$
being an exponential random variable independent of  Brownian motion
driving the diffusion $Y$ we find the closed formula of $\mathbb{E}
X_{2T_{\lambda}}^\alpha$ (Theorem \ref{tw:explicitmom}). In
Proposition \ref{prop:St} we give a closed formula for moments in
the Stein model.

Summing up, we present  forms  of some
interesting functionals of Brownian motion.
 Moreover, we find interesting links among a
GBM,  Markov processes arisen during generalization of the so called
Pitmann's $2M-X$ theorem (see \cite{Mat}, \cite{MatI}, \cite{MatII},
\cite{MatIII}) and a hyperbolic Bessel process (see \cite{RY}).
Finally, we compute the moments of the asset price process in the
lognormal stochastic volatility and Stein models.

\section{Properties of some functionals of geometric Brownian motion}
Let  $(\Omega,\mathcal{F},\mathbb{P})$ be a complete probability
space with filtration $\mathbb{F}= (\mathcal{F}_t)_{t \in
[0,\infty)}$ satisfying the usual conditions.
Let the process $Y$ be of the form
\begin{equation}  \label{def:Y}
 Y_t = \exp \Big(B_t -\frac{t}{2}\Big),
 \end{equation}
where  $B$ is a Brownian motion.
 Functionals of $Y$ play a crucial role in many
problems  of modern stochastic analysis. The studies of the
properties of integral $\int_0^tY_u^2du$ are  motivated by the
problem of pricing Asian options. (see \cite{Mans08}, \cite{DGY}).
The process $Y_t^{-2}\int_0^t Y^2_udu$ has been considered by
Matsumoto and Yor in several works concerning laws of Brownian
motion functionals. Along with $Y_t^{-1}\int_0^t Y^2_u du$ it plays
a central role in a generalization of the Pitmann's $2M-X$ theorem
(for details see for instance \cite{Mat}, \cite{MatI},
\cite{MatIII}).
 Here, we investigate, among others,  the properties of functional
$\Gamma$ defined, for $\beta>0$, by
\begin{equation}  \label{def:gamma}
 \Gamma_t = \frac{Y_t}{1+\beta\int_0^tY_sds}.
 \end{equation}
It turns out that this process plays a crucial role in the problem
of computing the moments of  the asset price in the lognormal stochastic volatility  model (see Section 3).
 We also find
some new properties of the exponential functional
 \begin{align}
  \label{def-A}
& A^{(\mu)}_t := \int_0^t(Y^{(\mu)}_u)^2du ,
\end{align}
where, for $\mu \in R$,
 \begin{align}
 \label{def-Y1} & Y^{(\mu)}_t := \exp \Big(B_t + \mu t\Big).
\end{align}
Therefore, $Y$ defined by \eqref{def:Y} is by definition equal to
$Y^{(-1/2)}$, so $\Gamma$ is a functional of  $Y^{(-1/2)}$. We also
consider the random variable (which is often called a perpetuity in
the mathematical finance literature):
\begin{align}
  \label{def-A-inf}
& A^{(\mu)}_\infty := \int_0^\infty (Y^{(\mu)}_u)^2du .
\end{align}
 We start from
investigation of $\Gamma.$

\subsection{Some properties of $\Gamma$}
\begin{proposition} \label{dyfuzja} If  $\ \Gamma$ is given by \eqref{def:gamma},
then $ \Gamma_0 = 1$ and
\begin{align} \label{eq:gamma}
d\Gamma_t = \Gamma_tdB_t - \beta \Gamma^2_t dt .
\end{align}
\end{proposition}
\begin{proof} It follows easily from the It\^o lemma.
\end{proof}
\begin{proposition} \label{recurence}
Let $\Gamma$ be given by \eqref{def:gamma} and $p_k(t) =
\int_0^t\mathbb{E}\Gamma_u^kdu$, $k\in \mathbb{Z}$,  $t \in [0, T]$.
Then  the sequence  of functionals $(p_k)$ satisfies the following
recurrences:
\begin{align}\label{momrec}
 p'_k(t) &= 1 + \frac{k(k-1)}{2}p_k(t) - \beta k p_{k+1}(t),
\end{align}
and
\begin{equation}  \label{def:p1}
 p_1(t) = \frac{1}{\beta} \mathbb{E}\Big(\ln(1 + \beta\int_0^tY_udu) \Big).
\end{equation}
\end{proposition}
\begin{proof}
By Proposition \ref{dyfuzja} and  the It\^o lemma  we have
\begin{align}\label{gamma}
\Gamma_t^k  = 1 + k \int_0^t\Gamma_u^kdB_u -
k\beta\int_0^t\Gamma_u^{k+1}du +
\frac{k(k-1)}{2}\int_0^t\Gamma^k_udu.
\end{align}
The local martingale $\int_0^t\Gamma_u^kdB_u$ is a true martingale
as
\begin{align*}
    \mathbb{E}\int_0^t\Gamma_u^{2k}du \leq \mathbb{E}\int_0^tY_u^{2k}du < \infty.
\end{align*}
Taking expectation of both sides of (\ref{gamma}) we obtain
\eqref{momrec}. Further
\begin{align*}
p'_1(t) &= \mathbb{E}\Big(\frac{Y_t}{1+\beta\int_0^tY_udu}\Big)
= \frac{1}{\beta}\mathbb{E}\frac{\partial}{\partial t}\Big(\ln \Big(1 + \beta\int_0^tY_udu\Big)\Big)\\
&= \frac{1}{\beta}\frac{\partial}{\partial t}\mathbb{E}\Big(\ln
\Big(1 + \beta\int_0^tY_udu\Big)\Big),
\end{align*}
as $\ln (1 + \beta\int_0^tY_udu) \leq \ln (1 + \beta\int_0^TY_udu)$
and $\mathbb{E}\ln (1 + \beta\int_0^TY_udu) <\infty$, which implies
\eqref{def:p1}.
\end{proof}
\begin{remark} \label{uw1}
Since, by \eqref{momrec},
\begin{align}\label{E-gamma}
\mathbb{E}\Gamma_t^k  = \frac{k(k-1)}{2}p_k(t) - \beta k p_{k+1}(t),
\end{align}
 Proposition \ref{recurence} allows to compute
$\mathbb{E}\Gamma_t^k$ for $k\in\mathbb{Z}$. Taking $k=-1$, we
easily obtain from  \eqref{momrec} that
\begin{align*}
    p'_{-1}(t) = 1+ \beta t + p_{-1}(t), \qquad
    p_{-1}(0) = 0.
\end{align*}
This solution is given by the formula $p_{-1}(t) = (\beta -1)e^t
+\beta t + 1 + \beta$. Notice that, having $p_{-1}$ we get
recursively from \eqref{momrec} the  functions $p_{-2}, p_{-3},...$.
Using the function $p_1$ we can establish  $p_2, p_3, ..$. So, using
\eqref{E-gamma}, we can find all moments  $\mathbb{E}\Gamma_t^k$ for
$k\in\mathbb{Z}$, provided we know the form of $p_1$. Therefore, to
finish this computation we need to find the closed form of the
function $p_1$. The function $p_1$ is given by \eqref{def:p1}, so we
have to find $\mathbb{E}\Big(\ln(1 + \beta\int_0^tY_udu) \Big).$ The
form of $p_1$ is presented in Corollary \ref{wn2}.
\end{remark}
Now we investigate the connection of $\Gamma$ with  hyperbolic
Bessel processes. Let us recall that a diffusion $R$ with the
generator given by
\begin{equation}
    \mathcal{A} = \frac{1}{2}\frac{d^2}{dx^2} + \big(\alpha+\frac{1}{2}\big)\coth(x)\frac{d}{dx},
\end{equation}
for $\alpha\in\mathbb{R}$, is called  a hyperbolic Bessel  (HB)
process with the parameter $\alpha$ (see \cite{RY} or \cite{B}).
Therefore $R$ satisfies
\begin{equation} \label{eq:hbp}
    dR_t = dB_t + \big(\alpha+\frac{1}{2}\big)\coth(R_t)dt.
\end{equation}
We express the Laplace'a transform  of functional  $\Gamma$ in terms
of the Laplace'a transform of $\cosh$ of $R$.
\begin{theorem} \label{tw:rephb} Let $R$ be  a
hyperbolic Bessel process with the parameter $\alpha= -1$ and
$\Gamma$ be given by \eqref{def:gamma}. For $\lambda\geq 0$ we have:
\begin{equation}
    \mathbb{E} e^{-\lambda \Gamma_t} = \mathbb{E}e^{-\beta(\cosh(R_t) - 1)},
\end{equation}
where the initial value of the process $R$ satisfies $\cosh(R_0) =
\frac{\lambda}{\beta} +1$.
\end{theorem}
\begin{proof}
Let $\theta_t=\beta\Gamma_t$. Then, by \eqref{eq:gamma}
\begin{equation}
    d\theta_t = \theta_tdZ_t -\theta_t^2dt
\end{equation}
and  $\theta_0 =  \beta $. Moreover, for $x\geq 0$
\begin{align*}
    de^{-x\theta_t} = -e^{-x\theta_t}(x\theta_tdZ_t - x\theta_t^2dt) +\frac{1}{2}e^{-x\theta_t}x^2\theta_t^2dt.
\end{align*}
Taking $p(t,x):= \mathbb{E}e^{-x\theta_t}$ we get from the last
expression that $p$ satisfies the PDE
\begin{equation}\label{eq:zag}
    \frac{\partial p}{\partial t} = \Big(x + \frac12x^2\Big)\frac{\partial^2 p}{\partial x^2},
\end{equation}
with $p(0,x) = e^{-x\beta}$. Therefore, the Laplace'a transform of
$\theta_t$ for $\lambda\geq 0$ is a solution of (\ref{eq:zag}).
Consider a stochastic differential equation (SDE)
\begin{equation} \label{difH}
    dH_t = \sqrt{H_t^2 + 2H_t}dB_t, \qquad H_0= \frac{\lambda}{\beta} \geq0.
\end{equation}
 Since for any $0\leq y \leq x$
\begin{align}
    \Big| \sqrt{x^2 + 2x} - \sqrt{y^2 +2y} \Big| \leq \sqrt{(x-y)^2 + 2(x-y)}
\end{align}
there exists a weak solution to SDE (\ref{difH}) and the trajectory
uniqueness holds for (\ref{difH}) (see \cite[Theorem 5.5.4]{KS} and
\cite[Theorem 5.40.1]{RW}).
 Thus, by the Feynman-Kac theorem (after changing terminal condition
to the initial one in the Cauchy problem  \eqref{eq:zag}) we obtain
that the function $u(t,x):= \mathbb{E}_xe^{-\beta H_t}$ is the
unique bounded solution of (\ref{eq:zag}) with $p(0,x) =
e^{-x\beta}$ (see
 \cite[Theorem 5.7.6]{KS}).
 Let us define the diffusion $S_t := H_t + 1$. It is easy to check
that
\begin{equation} \label{difS}
  dS_t =   \sqrt{S_t^2 -1 }\ dB_t, \qquad H_0= \frac{\lambda}{\beta} +1 .
\end{equation}
By the same arguments as before there exists a weak solution to
(\ref{difS}) and the trajectory uniqueness holds for (\ref{difS}).
 Now,  we
observe that the diffusion $U_t=\cosh(R_t)$, where $R$ is the
hyperbolic Bessel process  with the parameter $-1$, and such that
$\cosh(R_0) = a + 1$ is the solution of (\ref{difS}). So, the
processes $S$ and $U$ have the same law. Thus,
$$
 \mathbb{E} e^{-\lambda \Gamma_t} = \mathbb{E} e^{-\frac{\lambda}{\beta} \theta_t} =
p(t,{\lambda}/{\beta})= \mathbb{E} e^{-\beta H_t} = e^{\beta}
\mathbb{E} e^{-\beta S_t} = e^{\beta}
  \mathbb{E}e^{-\beta\cosh(R_t)}.
$$
 This ends the proof.
\end{proof}
The connections between hyperbolic Bessel processes and functionals
of geometric Brownian motion are presented in \cite{jw}.
\subsection{Some properties of $(1+\beta
A^{(\mu)}_t)$}

We start form the computation of the Laplace'a transform of
 $(1+\beta
A^{(\mu)}_t)^{-1}$.
 It is worth to remark that we compute it for a fixed time $t$.
We can find in literature (see for instance \cite{MatII}) that the
problem of computing of expectations for functionals of geometric
Brownian motion for a fixed time is in general much difficult than
with stochastic one (see also Subsection \ref{sec-ran}).

 Let us recall that a squared
$\delta$-dimensional radial Ornstein-Uhlenbeck process  with the
parameter $-\lambda$ for $\delta\geq0, \lambda\in\mathbb{R}$, is the
solution of the SDE
\begin{equation} \label{def-rOU}
    X_t = x+ \int_0^t(\delta-2\lambda X_s)ds + 2\int_0^t\sqrt{X_s}dW_s,
\end{equation}
where $W$ is a standard Brownian motion. For detailed studies of
these processes see \cite{AGY} and \cite{SB}. If $\lambda =0$, then
the strong solution of \eqref{def-rOU} is  a squared
$\delta$-dimensional Bessel process (see \cite{RY}). The number
$\delta/2 -1$ is called the index of the process.
 In the sequel we will use the notation $X^x$ for the process $X$
starting from $x$, i.e. $X_0=x$.

In the next two theorems we find a probabilistic representation of
the Laplace transform of $(1+\beta A^{(\mu)}_t)^{-1}.$
\begin{theorem} \label{tw:qu}
Assume $\beta \in(0,1]$,  $\mu \in \mathbb{R}$ and $ t > 0$.
Then, for any $\lambda >0$
\begin{equation} \label{eq:transf}
\mathbb{E}\exp\Big(-\frac{\lambda}{1+\beta A^{(\mu)}_t}\Big) =
\mathbb{E}\phi_t(\theta^{\lambda}(-\ln(\sqrt{\beta}))),
\end{equation}
where $\theta^{\lambda}(t)$ is a squared 0-dimensional radial
Ornstein-Uhlenbeck process with the parameter $-1$ such that
$\theta(0) = \lambda$ and, for $x >0$,
\begin{equation}
    \phi_t(x) = \psi_t(1,x),
\end{equation}
and for $x>0, s\geq 0$
\begin{align} \label{eq:postac-psi}
    \psi_t(s,x) = \mathbb{E} G_t(R^x(s/2)).
\end{align}
Here $R^x$ is a squared Bessel process with the index $-1$ starting
from $ x$, and
\begin{align} \label{eq:postac-G}
    G_t(x) &= e^{-t\mu^2/2}\mathbb{E}\exp\Big(\mu B_t +\frac{1}{2t}\Big(B_t^2 -
    \varphi^2_x(B_t)\Big)\Big),
 \end{align}
$B$ is a standard Brownian motion, and
\begin{align*}
    \varphi_x(y)
    &= \ln\Big(xe^{-y} + \cosh(y)+\sqrt{x^2e^{-2y}+\sinh^2(y) +2xe^{-y}\cosh(y)}\ \Big).
                            \end{align*}
\end{theorem}
\begin{proof}
Let $B$ be a standard Brownian motion under  $\mathbb{P}$. Let us
define the function
\begin{equation} \label{def-q}
q(s,x) := \mathbb{E}\exp\Big(-\frac{x}{1+s A^{(\mu)}_t}\Big)
 \end{equation}
for $x \geq 0$, $s \in [0,1]$. Observe that $q(s,x)\leq 1$. It is
not difficult to check, using the Lebesgue theorem, that $q$ belongs
to the class $C^{1,2}([0,\infty)\times [0,\infty))$. Moreover, it is
easy to see  that $q$ satisfies the partial differential equation
\begin{equation} \label{eq:dla-q}
    -s\frac{\partial q}{\partial s} =
x \Big(\frac{\partial q}{\partial x}+\frac{\partial^2 q}{\partial
x^2}\Big),
\end{equation}
and $q(s,0) = 1$. Define
$$\phi_t(x):= q(1,x) =
\mathbb{E}\exp\Big(-\frac{x}{1+ A^{(\mu)}_t}\Big) \quad {\rm and}
\quad
p(s,x):= q(e^{-s},x).$$
 Then $p(s,x)$ satisfies the partial
differential equation
\begin{equation} \label{eq:diffeqp}
    \frac{\partial p}{\partial s} = x \Big(\frac{\partial p}{\partial x}+\frac{\partial^2 p}{\partial x^2}\Big),
\end{equation}
 $s\geq 0, x\geq0$ and $p(0,x) = \phi_t(x)$. Consider  a diffusion $U$ with
 a generator of the form
\begin{equation}
    \mathcal{A}_U = x\frac{d^2}{dx^2} + x\frac{d}{dx}.
\end{equation}
This diffusion satisfies the SDE
\begin{align} \label{proc1}
    dU_{t} = \sqrt{2} \sqrt{U_{t}} dW_{t} + U_{t}dt,
\end{align}
where $W$ is a standard Brownian motion. Since  $p(s,x)\leq 1$ and
 $p \in C^{1,2}([0,\infty)\times [0,\infty))$, using the
Feynman-Kac theorem (after changing terminal condition to the
initial one in the Cauchy problem  \eqref{eq:diffeqp}) we obtain
that $p$ admits the stochastic representation $p(s,x) =
\mathbb{E}\phi_t(U^x_s)$  (see
 \cite[Theorem 5.7.6]{KS}).

 Observe that $\theta(s) :=
U_{2s}$ is a 0-dimensional radial Ornstein-Uhlenbeck process with
the parameter $-1$, as  \eqref{proc1} takes, by the scaling property
of Brownian motion, the form
\begin{equation}
    d\theta(t) = 2\sqrt{\theta(t)}dB_t + 2\theta(t)dt.
\end{equation}
Thus
\begin{align} \label{eq:27}
    q(s,x) = p(-\ln s,x) =
    \mathbb{E}\phi_t(\theta^x( - (1/2) \ln s))
    = \mathbb{E}\phi_t(\theta^x(-\ln\sqrt{s})),
\end{align}
 $\theta(0) = x$. Hence, taking $x= \lambda, s= \beta $ we see that
\eqref{eq:27} gives \eqref{eq:transf}:
$$
\mathbb{E}\exp\Big(-\frac{\lambda}{1+\beta A^{(\mu)}_t}\Big) =
q(s,\lambda) =
\mathbb{E}\phi_t(\theta^{\lambda}(-\ln(\sqrt{\beta}))).
$$
To finish the proof we have to find the form of  $\phi_t$. To do
this we define the new functions:
 \begin{align}  \label{def:postac-psi}
 \psi_t(s,x):= & \
\mathbb{E}\exp\Big(-\frac{x}{s + A^{(\mu)}_t}\Big), \\
\label{def:postac-G} G_t(x):= & \
\mathbb{E}\exp\Big(-\frac{x}{A^{(\mu)}_t}\Big),
\end{align}
for  $s\geq 0, x\geq0$. So $\psi_t(1,\lambda)=\phi_t(\lambda)$ and
$G_t(x)=  \psi_t(0,x)$. Observe that $\psi_t$ satisfies the partial
differential equation
\begin{equation} \label{eq:diffeqpsi}
    \frac{\partial\psi_t}{\partial s} = x\frac{\partial^2\psi_t}{\partial
    x^2}.
\end{equation}
 Consider  a diffusion $X$ with a generator of the form
\begin{align*}
    \mathcal{A}_X = x\frac{d^2}{dx^2}.
\end{align*}
Using again the Feynman-Kac theorem
we deduce that $\psi$ admits the stochastic
representation $\psi_t(s,x) = \mathbb{E}G_t(X^x_s)$.
 Now observe that $R(s) := X^x_{2s}$ satisfies
\begin{equation}
    dR(t) = 2\sqrt{R(t)}dB_t,
\end{equation}
so $R$ is a squared Bessel process with the index $-1$. Therefore,
we obtain \eqref{eq:postac-psi}.
It remains  to compute the form of the function $G_t$. Define a new
probability  measure $\mathbb{Q}$ by
\begin{equation}
    \frac{d\mathbb{Q}}{d\mathbb{P}}\Big|_{\mathcal{F}_t} = \exp\Big(-\mu B_t -\frac{\mu^2}{2}t\Big).
\end{equation}
Since $B$ is a standard Brownian motion under $\mathbb{P}$, then
$\widehat{B}_t = B_t+\mu t$ is a standard Brownian motion under
$\mathbb{Q}$, by
 the Girsanov theorem. For ${\widehat A}^{(0)}_t := \int_0^t
e^{2\widehat{B}_u}du$
 we have $A^{(\mu)}_t={\widehat A}^{(0)}_t $, so
\begin{align*}
    G_t(x) &= \mathbb{E}e^{-\frac{x}{A^{(\mu)}_t}}= \mathbb{E}_{\mathbb{Q}}
    \Big(e^{-\frac{x}{A^{(\mu)}_t}}e^{\mu B_t
    +\frac{\mu^2}{2}t} \Big) =
    e^{-\mu^2 t/2}\mathbb{E}_{\mathbb{Q}}\exp\Big(\mu\widehat{B}_t -\frac{x}{{\widehat A}^{(0)}_t}\Big).
\end{align*}
 Now we use the   Matsumoto-Yor result \cite[Thm. 5.6]{MatII},
which states that
\begin{align*}
    \mathbb{E}_{\mathbb{Q}}\Big(\exp\Big(-\frac{x}{{\widehat A}^{(0)}_t}\Big)\Big|
     \widehat{B}_t = y\Big) = \exp\Big(-\frac{\varphi^2_x(y)-y^2}{2t}\Big),
\end{align*}
where
\begin{align*}
    \varphi_x(y) &= \arg\cosh(xe^{-y} + \cosh(y))\notag\\
                            &= \ln\Big(xe^{-y} + \cosh(y)+\sqrt{x^2e^{-2y}+\sinh^2(y) +2xe^{-y}\cosh(y)}\Big).
    \end{align*}
In result we obtain
\begin{align*}
    G_t(x) &= e^{-t\mu^2/2}\mathbb{E}\exp\Big(\mu B_t +
    \frac{1}{2t}\Big(B_t^2 - \varphi^2_x(B_t)\Big)\Big),
\end{align*}
where $B$ is a standard Brownian motion under $\mathbb{P}$. This
finishes the proof.
\end{proof}
Our next theorem provides another probabilistic representation
of Laplace transform of $(1 +\beta A_t^{(\mu)})^{-1}$, now for
$\beta>0$.
\begin{theorem} \label{tw:freeb}
Fix $\beta>0, \mu \in\mathbb{R}$ and $t\geq 0$. Then, for any
$\lambda \geq 0$,
\begin{equation} \label{generalb}
\mathbb{E}\exp\Big(-\frac{\lambda}{1+\beta A^{(\mu)}_t}\Big) =
\mathbb{E}G_t(R^{\lambda/\beta}(1/(2\beta) ),
\end{equation}
where $R^{\lambda/\beta}$ is a squared Bessel process with the index
${\lambda/\beta}$ starting from $ x$ and $G_t$ is defined by
\eqref{eq:postac-G}.
\end{theorem}
\begin{proof} Formula (\ref{generalb}) follows from
the proof of Theorem \ref{tw:qu}, since for $s\geq 0, x \geq 0$ we
have, by \eqref{def:postac-psi} and \eqref{eq:postac-psi},
\begin{equation}
    \mathbb{E}\exp\Big(-\frac{x}{s + A_t^{(\mu)}}\Big) = \psi_t (s,x) = \mathbb{E}G_t(R^x(s/2))
\end{equation}
and taking $x = \lambda / \beta$ and $s = 1/ \beta$ we obtain
(\ref{generalb}).
\end{proof}
\begin{corollary}
 For $\beta\in(0,1]$ we have
\begin{equation}\label{Rousb}
\mathbb{E}\phi_t(\theta^{\lambda}(-\ln(\sqrt{\beta}))) =
\mathbb{E}G_t(R^{\lambda/\beta}(1/(2\beta) ),
\end{equation}
where $\theta^{\lambda}$ is defined in Theorem \ref{tw:qu}.
\end{corollary}
\begin{proof}
 (\ref{Rousb}) follows from (\ref{generalb}) and
\eqref{eq:transf}.
\end{proof}

Using this result we can obtain the expectations of $\ln (1+\beta
A^{(\mu)}_t)$ and $(1+\beta A^{(\mu)}_t)^{-1}$.
\begin{theorem} \label{1st} Fix $\beta >0$,  $\mu \in \mathbb{R}$ and $t\geq 0$. Then
\begin{align} \label{eq:nalog}
    \mathbb{E}\ln (1+\beta A^{(\mu)}_t) & =
    \int_0^{\infty}\frac{G_t(y)}{y}(1-e^{-y\beta})dy, \\
    \label{psiinzerooo}
\mathbb{E}\Big(\frac{1}{1 +\beta A_t^{(\mu)}} \Big) & = 1 -
\beta\int_0^{\infty}G_t(y)e^{-y\beta}dy,
\end{align}
where $G_t$ is given by \eqref{eq:postac-G}.
\end{theorem}
\begin{proof} Let $f(\beta) = \mathbb{E}\ln (1+\beta A^{(\mu)}_t)$ for $\beta>0$.  Since
\begin{align*}
    \mathbb{E}|\ln (1+\beta A^{(\mu)}_t)| \leq 1 + \beta\mathbb{E}A^{(\mu)}_t <
    \infty,
\end{align*}
the  function $f$  is  well defined, continuous and $f(0)=0$.
 Moreover, for $\beta>0$,
\begin{align}\label{fprim}
   f'(\beta) =
 \mathbb{E}\Big(\frac{A^{(\mu)}_t}{1+\beta A^{(\mu)}_t}\Big) =
 \frac{1}{\beta}\Big(1 - \mathbb{E}\Big(\frac{1}{1+\beta A^{(\mu)}_t}\Big)\Big).
\end{align}
By definition of $\psi$ (see \eqref{def:postac-psi})  and
\eqref{eq:postac-psi} we know that for $s\geq 0,$ $ x \geq 0$ we
have
\begin{equation} \label{psiR}
    \mathbb{E}\exp\Big(-\frac{x}{s + A_t^{(\mu)}}\Big) = \
     \psi_t (s,x) = \mathbb{E}G_t(R^x(s/2)),
\end{equation}
where $R^x$ is a Bessel process with the index $-1$ starting from
$x$. Since, by definition (see \eqref{def:postac-G}), $G_t(0) = 1$
and the transition density functions for the process $R^x$ are known
(see \cite[Chapter IX, Corollary 1.4]{RY}) we can write
\begin{align} \label{TFG}
    \mathbb{E}G_t(R^x(s/2)) &= e^{-x/ s} + \int_0^{\infty}G_t(y)\frac{1}{s} \sqrt{\frac{x}{y}} e^{-(x+y)/s}I_1(2\sqrt{xy}/s)dy,
\end{align}
where $I_1$ is the modified Bessel function. Let us recall that
$(I_1(x)/x)' = I_2(x)/x$ (see \cite[Appendix 2]{SB}). Hence and by
(\ref{TFG}) and (\ref{psiR}) we obtain, for $x \in [0,1]$,
\begin{align*}
 - \mathbb{E}&\Big(\frac{1}{s + A_t^{(\mu)}}\exp\Big(-\frac{x}{s + A_t^{(\mu)}}\Big) \Big)
 = \ \frac{\partial\psi(s,x)}{\partial x} \\
&= -\frac{1}{s}e^{-x/s} + \frac{\partial}{\partial x}
\Big(\int_0^{\infty}G_t(y)\frac{1}{s} \sqrt{\frac{x}{y}} e^{-(x+y)/s}I_1(2\sqrt{xy}/s)dy\Big)\\
&= -\frac{1}{s}e^{-x/s} + \frac{2e^{-x/s}}{s^2}(1 - \frac{x}{s})\Big(\int_0^{\infty}G_t(y)
e^{-y/s}(2\sqrt{xy}/s)^{-1}I_1(2\sqrt{xy}/s)dy\Big)\\
&+ \frac{2e^{-x/s}}{s^2}\Big(\int_0^{\infty} \sqrt{x}\ G_t(y)
e^{-y/s}(2\sqrt{xy}/s)^{-1}I_2(2\sqrt{xy}/s)
\frac{\sqrt{y}}{s}dy\Big).
\end{align*}
Let  $x$ tend to $0$. Since we can pass with the limit under the
integrals, using the asymptotic behavior of Bessel functions, i.e.
$I_1(x)/x \backsimeq 1/2$ and $I_2(x)/x \backsimeq x/8$ in
neighborhood of $0$  (see \cite[Appendix 2]{SB}) we obtain
\begin{align} \label{psiinzero}
\mathbb{E}\Big(\frac{1}{s + A_t^{(\mu)}}\Big) = \frac{1}{s} -
\frac{1}{s^2}\int_0^{\infty}G_t(y)e^{-y/s}dy,
\end{align}
or equivalently
\begin{align} \label{psiinzeroo}
\mathbb{E}\Big(\frac{1}{1 +\frac{1}{s} A_t^{(\mu)}} \Big) = 1 -
\frac{1}{s}\int_0^{\infty}G_t(y)e^{-y/s}dy.
\end{align}
Putting $s = 1/\beta$ in (\ref{psiinzeroo}) yields
\eqref{psiinzerooo}. From (\ref{fprim}) and (\ref{psiinzerooo}) we
conclude
\begin{align*}
    f'(\beta) = \int_0^{\infty}G_t(y)e^{-y\beta}dy.
\end{align*}
This finishes proof of the theorem,  since $f(0) = 0$.
\end{proof}
\begin{remark} \label{uw2a}
 Formula \eqref{psiinzerooo} gives
the closed expression of $\mathbb{E}((1 +\beta A_t^{(\mu)})^{-1})$
for $\beta>0$. The density of $A_t^{(\mu)}$ is known in literature,
but due to complicated nature of Hartman-Watson distribution, it can
hardly be used for numerical computations (see for instance
\cite{MatII} and \cite{BRY}). Since the simple form of function
$G_t$ is   given explicitly, the formulae (\ref{eq:nalog}) and
(\ref{psiinzerooo}) allows to obtain numerically $ \mathbb{E}\ln
(1+\beta A^{(\mu)}_t)$ and $\mathbb{E}((1 +\beta
A_t^{(\mu)})^{-1})$.
\end{remark}
 Theorem \ref{1st} allows to find the first function $p_1(\cdot)$
for the recurrence established in Proposition \ref{recurence}.
\begin{corollary} \label{wn2}
Let $p_1$ be given by \eqref{def:p1}. Then
\begin{align}
 p_1(t) = \frac{1}{\beta} \int_0^{\infty}\frac{G_{t/4}(y)}{y}(1-e^{-4\beta y})dy,
\end{align}
where $G$ is defined by \eqref{eq:postac-G}.
\end{corollary}
\begin{proof}
Since $Z_u=(1/2)B_{4u}$ is  a standard Brownian motion, we infer
that
\begin{align*}
  p_1(4t) &= \mathbb{E}\ln(1 + \beta\int_0^{4t}Y_udu) =
  \mathbb{E}\ln(1 + 4\beta\int_0^{t}e^{B_{4u} - 2u}du)\notag\\
    &= \mathbb{E}\ln(1 + 4\beta\int_0^{t}e^{2(Z_u - u)}du)
    = \mathbb{E}\ln (1+ 4\beta A^{(-1)}_t).
\end{align*}
\end{proof}
Let us now observe that we can deduce more general fact from the
 proof of
Theorems \ref{tw:qu} and \ref{1st}. It turns out
 that for an arbitrary strictly positive random variable we can find a
representation of the Laplace'a transform of $(s+\xi)^{-1}$ for
$s\geq 0$, in terms of a  squared Bessel process.
 It gives also a second simple proof of general version of
\eqref{eq:postac-psi}.
\begin{theorem}
\label{tw:og-rep} Let $\xi$ be a  strictly positive random variable.
Then for any $x\geq 0, s \geq 0$
\begin{equation}
    \mathbb{E}\exp\Big(-\frac{x}{s + \xi}\Big) = \mathbb{E}G(R^x(s/2)),
\end{equation}
where $R^x$ is a squared Bessel process with the index $-1$ starting
from $x$ and
\begin{equation} \label{general-G}
G(x) = \mathbb{E}\exp\Big(-\frac{x}{\xi}\Big).
\end{equation}
Moreover for $\beta \geq 0$
\begin{equation} \label{generallog}
    \mathbb{E}\ln (1+\beta \xi) = \int_0^{\infty}\frac{G(y)}{y}(1-e^{-y\beta})dy.
\end{equation}
\end{theorem}
\begin{proof} Let's take a copy of $R^x$ independent of $\xi$. Then
\eqref{general-G} implies
\begin{align*}
     \mathbb{E}G(R^x(s/2)) &=  \mathbb{E}\exp\Big(-\frac{R^x(s/2)}{\xi}\Big)
     = \mathbb{E}\mathbb{E}\Big(\exp\Big(-\xi^{-1}R^x(s/2)\Big)\Big|\xi\Big)\\
     &= \mathbb{E}\exp\Big(-\frac{x\xi^{-1}}{1+\xi^{-1}s}\Big)
     = \mathbb{E}\exp\Big(-\frac{x}{s + \xi}\Big),
\end{align*}
where we used the form of Laplace transform of squared Bessel
process (see \cite[Chapter XI, page 441]{RY})). The proof of
(\ref{generallog}) goes in the same way as in Theorem \ref{1st}.
\end{proof}
Now, we use formula \eqref{psiinzerooo} and the results of Matsumoto
and Yor to obtain some interesting connections between
$\mathbb{E}((1 +\beta A_t^{(\mu)})^{-1})$ and the conditional
expectation of functionals of geometric Brownian motion with
opposite drift.
\begin{proposition} \label{MMA} For $\mu >0$ and $\beta>0$ we have
\begin{equation}
    \mathbb{E}\Big(\frac{1}{1+2\beta A_t^{(\mu)}} \Big) = 1 - 2\beta\mathbb{E}\Big(A_t^{(-\mu)}|
    A_{\infty}^{(-\mu)} = 1/(2\beta)\Big) .
\end{equation}
\end{proposition}
\begin{proof} By the result of Matsumoto and Yor \cite[Thm. 2.2]{Mat-I}
the process $\{B_t^{(-\mu)}, t\geq 0\}$ on the set
$\{A_{\infty}^{(-\mu)} = 1/(2\beta)\}$ has the same distribution as
the process $\{B_t^{(\mu)} - \log (1 + 2\beta A_t^{(\mu)}), t\geq
0\}$ for $\mu > 0$.
From that we obtain
\begin{align*}
    \mathbb{E}\Big(A_t^{(-\mu)}|A_{\infty}^{(-\mu)} = 1/(2\beta)\Big) =
    \mathbb{E}\int_0^t\frac{e^{2B_s^{(\mu)}}}{(1+2\beta A_s^{(\mu)})^2}\ ds
    = \frac{1}{2\beta}\Big(1-\mathbb{E}\Big(\frac{1}{1+2\beta A_t^{(\mu)}}\Big)\Big).
\end{align*}
\end{proof}
\begin{proposition} Let $\beta>0$
and $\mu\in\mathbb{R}$. Then
\begin{equation}
\mathbb{E}\Big(A_t^{(-\mu)}| A_{\infty}^{(-\mu)} = 1/(2\beta)\Big) =
\frac{1}{2}\int_0^{\infty}G_t(y)e^{-y\beta}dy.
\end{equation}
\end{proposition}
\begin{proof} It follows from Proposition \ref{MMA} and \eqref{psiinzerooo}.
\end{proof}
\begin{proposition} \label{FWG} For $\beta>0, \mu > 0$ we have
\begin{equation} \label{eq:FWG}
            \mathbb{E} \Big( \frac{e^{2\mu B_t^{(-\mu)}}}{1+2\beta A_t^{(-\mu)}} \Big)
            = 1 - 2\beta \mathbb{E}\Big(A_t^{(-\mu)}|A_{\infty}^{(-\mu)} = 1/(2\beta)\Big) .
\end{equation}
\end{proposition}
\begin{proof} Fix $\mu >0$. Define the new probability measure $\mathbb{Q}$ by
\begin{equation}
    \frac{d\mathbb{Q}}{d\mathbb{P}}\Big|_{\mathcal{F}_t} = e^{-2\mu B_t - 2\mu^2 t}.
\end{equation}
The process $V_t = B_t +2\mu t$ is a standard Brownian motion under
$\mathbb{Q}$, so
\begin{align*}
    \mathbb{E} \Big( \frac{1}{1+2\beta A_t^{(\mu)}}  \Big) = \mathbb{E}_{\mathbb{Q}} \Big(
    \frac{e^{ 2\mu(V_t-\mu t) }}{1+2\beta\int_0^te^{2(V_u-\mu u)}du}  \Big)
    = \mathbb{E} \Big( \frac{e^{2\mu B_t^{(-\mu)}}}{1+2\beta A_t^{(-\mu)}} \Big),
\end{align*}
 Now the thesis follows from  Theorem \ref{MMA}.
\end{proof}
\begin{proposition} \label{GMM} For $t\geq 0$ we have
\begin{equation}
    p_1(t) = t - 4\beta\int_0^t\mathbb{E}\Big(A_{s/4}^{(-1)}|A_{\infty}^{(-1)} = 1/(4\beta) \Big)ds.
\end{equation}
\end{proposition}
\begin{proof}
We have $p'_1(t) = \mathbb{E}\Gamma_t$, $p_1(0)=0$, and
\begin{equation}
    \Gamma_{4t} = \frac{e^{2\overline{B}^{(-1)}_t}}{1+4\beta \overline{A}_t^{(-1)}},
\end{equation}
where $\overline{B}_t = B_{4t}/2$ is a standard Brownian motion, and
$\overline{A}^{(-1)}$ is defined by \eqref{def-A} with
$\overline{B}$ instead of $B$. Since
$$
\mathbb{E}\Big(A_{t/4}^{(-1)}|A_{\infty}^{(-1)} = 1/(4\beta) \Big) =
\mathbb{E}\Big(\overline{A}_{t/4}^{(-1)}|
\overline{A}_{\infty}^{(-1)} = 1/(4\beta) \Big),
$$
 the thesis follows from \eqref{eq:FWG} with $\mu=1$.
\end{proof}

\begin{remark} Notice, that we establish all the results for fixed
$t$. In  several papers (for instance \cite{Mans08}, \cite{MatII})
the integral functionals of a geometric Brownian motion  with random
time given by  random variable independent of Brownian motion and
with exponential distribution were investigated.
  In particular,
\begin{equation}
\mathbb{E}\ln\Big(1 +\beta\int_0^{T_{\lambda}}Y_u^2du\Big) =
\mathbb{E}\ln\Big(1+\beta\frac{\zeta_{1,a}}{\gamma_b}\Big),
\end{equation}
because
$\int_0^{T_{\lambda}}Y_u^2du =^d\frac{\zeta_{1,a}}{\gamma_b}$, where
$\zeta_{1,a}$ is a random variable with  beta distribution with the
parameters $1$ and $a = \frac{\sqrt{2\lambda +1/4}-1/2}{2}$,
 $\gamma_b$ is a random variable with
 gamma  distribution with the parameter $b = \frac{\sqrt{2\lambda +1/4}
+1/2}{2}$, $\zeta_{1,a}$ and $\gamma_b$ are independent (see
\cite{MatII}). Later, we also explore idea of using random time. In
subsection \ref{sec-ran} we show how to compute the moments in a
lognormal stochastic volatility model with random time being
exponentially distributed and independent of Brownian motion driving
the model.
\end{remark}

\section{Moments of the asset price in the lognormal stochastic
volatility and Stein models }

\subsection{Model of market}

We consider a market defined on a complete probability space
$(\Omega,\mathcal{F},\mathbb{P})$ with filtration $\mathbb{F}=
(\mathcal{F}_t)_{t \in [0,T]}$, $T<\infty$, satisfying the usual
conditions and  $\mathcal{F}=\mathcal{F}_T$. Without loss of
generality we assume the savings account to be constant and
identically equal to one. Moreover, we assume that  the price $X_t$
of the underlying asset at time $t$ has a stochastic volatility
$Y_t$ being a geometric Brownian motion  or  an Ornstein-Uhlenbeck
process, so the dynamics of the proces $X$ is given by
\begin{align}
    dX_t &= Y_tX_tdW_t, \label{defX},
\end{align}
where $X_0 = 1$.
 In case of $Y$ being a GBM  the dynamics of the vector $(X, Y)$
is given by \eqref{defX} and
\begin{align}
    dY_t &= Y_t dZ_t, \qquad Y_0 = 1 \label{defYc},
\end{align}
and in  case of $Y$ being an OU  the dynamics of the vector $(X, Y)$
is given by \eqref{defX} and
\begin{align}
    dY_t &= - \lambda Y_t dt + dZ_t,  \qquad Y_0 = 1 \label{defY-2}
\end{align}
for $\lambda > 0$. The processes $W,Z$ are correlated Brownian
motions, $d{\left\langle W,Z\right\rangle}_t = \rho dt$ with
$\rho\in[-1,1]$. In the both cases the process $X$ has the form
\begin{equation}
  X_t = e^{\int_0^tY_udW_u - \zfrac{\int_0^tY_u^2du}{2}} , \label{postacX}
\end{equation}
and this is a unique strong solution of SDE \eqref{defX} on $[0,T]$.
The existence and uniqueness follow directly from the assumptions on
$Y$ and the well known properties of stochastic exponent (see, e.g.,
Revuz and Yor \cite{RY}). Since the process $X$ is a local
martingale, there is no arbitrage on the market so defined.\\
 Notice that we can represent $W$ as
  \begin{equation} \label{eq:repr-W}
  W_t
= \rho Z_t + \sqrt{1-\rho^2}V_t,
\end{equation}
where $(V,Z)$ is a standard two-dimensional Wiener process. Using
(\ref{postacX}) and  \eqref{eq:repr-W} we can expressed the moment
of order $\alpha$ of  $X$ as
\begin{align} \label{potega}
& \mathbb{E}X_t^{\alpha} = \mathbb{E}e^{\alpha\rho\int_0^tY_udZ_u +
\alpha\sqrt{1-\rho^2}\int_0^tY_udV_u
-\frac{\alpha}{2}\int_0^tY_u^2du}\\ \notag & =
\mathbb{E}e^{\alpha\rho\int_0^tY_udZ_u +
\frac{\alpha^2(1-\rho^2)-\alpha}{2}\int_0^tY_u^2du} .
\end{align}

\subsection{Moments  of the asset price in the lognormal stochastic
volatility  model }

 \subsubsection{{\bf Moments of order $\alpha>0$}}

In this subsection we will calculate moments of order $\alpha>0$ in
the lognormal stochastic volatility model, so for $Y$ of the form
\eqref{defYc}. Jourdain \cite{Jou} gave a sufficient condition on
$\alpha > 1$ for existing of moments of $X$. Namely, Jourdain proved
that for $\alpha \in (1, (1 - \rho^2)^{-1})$ and $\rho \neq 0$ the
moments $EX^\alpha$ exist, but he didn't find the value of these
moments. We calculate  the value of moments for  $\alpha > 0$,
$\alpha (1 - \rho^2) < 1$ and $\rho\in[-1,1]$. Sin \cite{SIN}
established that the process $X$ is a true martingale if and only if
$\rho\leq 0$. First, we prove that the moment of order $\alpha$ of
the strong solution of \eqref{defX} is equal to the Laplace'a
transform of  the process $\Gamma$.
\begin{theorem}\label{tw:momrep1}
Let $t\in [0,T]$, $\alpha>0$, $\alpha (1 - \rho^2) < 1$
and  $\Gamma$ be given by \eqref{def:gamma}. If $X$ is given by
\eqref{defX}, then
\begin{align}
\label{eq:momrep} \mathbb{E}X_t^{\alpha} =
e^{-(\beta+\rho\alpha)}\mathbb{E}\exp\Big((\beta+\rho\alpha)\Gamma_t
\Big),
   \end{align}
where
\begin{equation}  \label{def:beta}
 \beta = \sqrt{\alpha - \alpha^2(1-\rho^2)}.
\end{equation}
\end{theorem}
\begin{proof}
Define a measure $\mathbb{Q}$  by
\begin{equation}
    \frac{d\mathbb{Q}}{d\mathbb{P}}\Big|_{\mathcal{F}_T} =
    e^{-\beta\int_0^TY_udZ_u -\frac{\beta^2}{2}\int_0^TY_u^2du},
\end{equation}
where $\beta$ is given by \eqref{def:beta}. The measure $\mathbb{Q}$
is a probability measure since, by \eqref{defYc},
\begin{align*}
    e^{-\beta\int_0^TY_udZ_u -\frac{\beta^2}{2}\int_0^TY_u^2du}&=
    e^{-\beta(Y_T-1) -\frac{\beta^2}{2}\int_0^TY_u^2du}\leq e^{\beta}.
\end{align*}
 Using (\ref{potega}) and the definition of $\mathbb{Q}$ we infer
\begin{align} \label{potega-GBM}
 \mathbb{E}X_t^{\alpha} = \mathbb{E}e^{\alpha\rho\int_0^tY_udZ_u +
\frac{\alpha^2(1-\rho^2)-\alpha}{2}\int_0^tY_u^2du} =
\mathbb{E}_{\mathbb{Q}}e^{(\rho\alpha + \beta)(Y_t - 1)}
.
\end{align}
By the Girsanov theorem, $B_t=Z_t + \int_0^t \beta Y_s ds$ is a
standard Brownian motion under $\mathbb{Q}$ and
\begin{align}\label{uklad}
    dZ_t = dB_t - \beta Y_t dt, \qquad
    Z_0 = 0.
\end{align}
 We know,  by the result of Alili,
Matsumoto and  Shiraishi \cite[Lemma 3.1]{AMY}, that the unique
strong solution of (\ref{uklad}) is given by
\begin{align*}
    Z_t = \frac{t}{2} +
    \ln\Big(\frac{U_t}{1+\beta\int_0^tU_sds}\Big),
    \end{align*}
    where
    $$
    U_t = e^{B_t -\frac{t}{2}}.
    $$
Therefore
\[
Y_t = \exp\Big(Z_t-\frac{t}{2}\Big) =
\frac{U_t}{1+\beta\int_0^tU_sds}.
\]
 The law of the process
$Y$ under $\mathbb{Q}$ is equal to the law of the process $\Gamma$
under $\mathbb{P}$, since the law of the process $Y$ under
$\mathbb{P}$ is equal to the law of the process $U$ under
$\mathbb{Q}$.
 Hence and  by \eqref{potega-GBM} we obtain  \eqref{eq:momrep}.
 This ends the proof.
\end{proof}
\begin{remark}
a) From Theorem \ref{tw:momrep1} we immediately see that all moments
exist provided $\rho^2=1$. \\
b) The condition $\alpha (1 - \rho^2) < 1$ is not a necessary
condition for existence of moments since in case of $\rho=0$ the
process $X$ is a martingale, so $EX_t$ exists. Although, for
$\rho=0$, $\mathbb{E}X_t^{\alpha}=\infty$ for $\alpha>1$
(\cite{Jou}).
\end{remark}
\begin{remark}
In Theorem \ref{tw:momrep1} we prove that the computation of
$\mathbb{E}X_t^{\alpha}$ for $\alpha>0$ such that $\alpha(1 -
\rho^2) < 1$ is equivalent to the computation of the  Laplace'a
transform of $\Gamma_t$ at point $\lambda=\beta+\rho\alpha$. In
turn,  the recurrence from Proposition \ref{recurence} allows to
find $\mathbb{E}\Gamma_t^k$, so we can find an approximation of the
Laplace'a transform of $\Gamma$ in the neighborhood of zero by its
moments, namely $\mathbb{E}e^{\lambda \Gamma_t}\approx \sum_{i=0}^N
\lambda^i\frac{\mathbb{E}\Gamma_t^i}{i!}$ for sufficiently large
$N$. In this way we obtain an approximate value of
$\mathbb{E}X_t^{\alpha}$.
 \end{remark}
 We can also express moments of
order  $\alpha>1$ in terms of the hyperbolic Bessel process with the
parameter $-1$.
\begin{theorem} \label{tw:repr-Bess}
Assume that $\alpha>1$,
$\alpha(1 - \rho^2) < 1$. Let $X$ be given by \eqref{defX} with $Y$
given by \eqref{def:Y} and  $R$ be a hyperbolic Bessel process with
the parameter $-1$. Then
\begin{equation}
    \mathbb{E}X_t^{\alpha} = e^{-\rho\alpha}\mathbb{E}e^{-\beta\cosh(R_t)},
\end{equation}
where $\cosh(R_0) = -\rho\alpha / \beta$, $\beta = \sqrt{\alpha -
\alpha^2(1-\rho^2)}$.
\end{theorem}
\begin{proof}
We use Theorem \ref{tw:momrep1} and Theorem \ref{tw:rephb} with
$\lambda = -(\beta +\rho\alpha)$,
and $\lambda>0$ provided $\alpha>1$. \end{proof}
\subsubsection{{\bf Moments with independent random time}} \label{sec-ran}

In this subsection we find the closed formulae for the moments in a
lognormal stochastic volatility, when the time is an exponential
random variable independent of  Brownian motion driving the
diffusion $Y$. The idea of considering such a time is not new and
can be find in many studies of Asian options (see for instance
\cite{Mans08}, \cite{MatII}).
\begin{proposition} \label{ConP}
Let   $T_{\lambda}$ be a random variable with exponential
distribution with the parameter $\lambda>0$. Assume that
$T_{\lambda}$ is independent of  a standard Brownian motion $B$.
 Let $Z_t
= 2B_{\frac{t}{4}}$, $Y_t = e^{-\frac{t}{2} + Z_t}$ and     $U_t =
e^{B_t-t}$. Then
\begin{align} \label{firstmom}
    \mathbb{E}&\Big(\ln\Big(1+\beta\int_0^{4T_{\lambda}}Y_udu\Big)\Big)
    =\frac{4 \beta}{\lambda}-4\beta^2\int_0^{\infty}
    \mathbb{E}(\int_0^{T_{\lambda}}U_s^2ds-K/4)^+(1+\beta K)^{-2}dK, \\
    \label{mom-U}
    \mathbb{E}&(\int_0^{T_{\lambda}}U_s^2ds-K/4)^+ = \frac{1}
{\lambda
\Gamma\Big(\frac{\sqrt{2\lambda+1}-1}{2}\Big)}\int_0^{2/K}e^{-u }
u^{\frac{\sqrt{2\lambda+1}-3}{2}}(1-Ku/2)^{\frac{\sqrt{2\lambda+1}+1}{2}}du.
\end{align}
\end{proposition}
\begin{proof}
It is obvious that
 $Y_{4t} = e^{-2t + Z_{4t}} =  e^{-2t + 2B_t} = U_t^2$.
Using the Taylor theorem with integral remainder to the
 function $f(x) = \ln(1 + \beta x)$ gives
\begin{align}
    \ln(1 + \beta x) = \beta x -\beta^2\int_0^{\infty}(x-K)^+(1+\beta
    K)^{-2}dK .
\end{align}
Hence replacing $x$ by  $\int_0^{4T_{\lambda}}Y_udu$ and taking
expectation we get
\begin{align} \label{eq1}
    \mathbb{E}&\Big(\ln(1 + \beta \int_0^{4T_{\lambda}}Y_udu)\Big) =
\beta\mathbb{E}\Big(\int_0^{4T_{\lambda}}Y_udu \Big)\\ \notag &
-\beta^2\int_0^{\infty}
\mathbb{E}(\int_0^{4T_{\lambda}}Y_udu-K)^+(1+\beta K)^{-2}dK \\
\notag
    & =  4\beta\mathbb{E}\Big(\int_0^{T_{\lambda}}U^2_sds\Big)
    -4\beta^2\int_0^{\infty}\mathbb{E}(\int_0^{T_{\lambda}}U_s^2ds-K/4)^+(1+\beta K)^{-2}dK.
\end{align}
Let $A_t = \int_0^tU_s^2ds$. The Mansuy and Yor theorem \cite[Thm.
6.1]{Mans08}  gives \eqref{mom-U} and
$\mathbb{E}A_{T_{\lambda}} = 1/\lambda$.
This and \eqref{eq1} completes the proof.
\end{proof}
In the next theorem we establish the explicit formula for moments of
$X_{T_{2\lambda}}$.
\begin{theorem}\label{tw:explicitmom}
Let $\alpha>0$, $\alpha (1 - \rho^2) < 1$ and  $T_{\lambda}$ be a
random variable with exponential distribution with the parameter
$\lambda>0$. Assume that $T_{\lambda}$ is independent of  Brownian
motions $V$ and $Z$ driving the process $X$. Then
\begin{align}
    \mathbb{E}&X_{2T_{\lambda}}^{\alpha}
    = \frac{1}{\lambda}
    e^{-(\alpha\rho +\beta)}\frac{\Gamma((1+\sqrt{4\lambda + 1})/2)}{\Gamma(1+\sqrt{4\lambda + 1})}\times \\
    \notag \Big(&\phi_1(1/2\beta)\int_0^{1/2\beta}
    e^{\frac{\alpha\rho-\beta}{2\beta}\frac{1}{y}}\phi_2(y)dy+
     \phi_2(1/2\beta)\int_{1/2\beta}^{\infty}e^{\frac{\alpha\rho-\beta}{2\beta}\frac{1}{y}}\phi_1(y)dy\Big),
\end{align}
where $\beta = \sqrt{\alpha - \alpha^2(1-\rho^2)}$,
\begin{align*}
    \phi_1(x) = x^{-(1+\sqrt{1+4\lambda})/2}\Phi\Big((1+\sqrt{1+4\lambda})/2,1+\sqrt{1+4\lambda},x^{-1}\Big),\\
    \phi_2(x) = x^{-(1+\sqrt{1+4\lambda})/2}\Psi\Big((1+\sqrt{1+4\lambda})/2,1+\sqrt{1+4\lambda},x^{-1}\Big),
\end{align*}
and $\Phi,\Psi$ denote the confluent hypergeometric functions of the
first and second kind, respectively
\begin{align*}
    \Phi(\alpha,\gamma,z) &= \sum_{k=0}^{\infty}\frac{(\alpha)_k}{(\gamma)_k}\frac{z^k}{k!} ,\\
    \Psi(\alpha,\gamma,z) &= \frac{\Gamma(1-\gamma)}{\Gamma(1+\alpha - \gamma)}\Phi(\alpha,\gamma,z)
     + \frac{\Gamma(\gamma - 1)}{\Gamma(\alpha)}z^{1-\gamma}\Phi(1+\alpha-\gamma,2-\gamma,z),
\end{align*}
where $(\alpha)_0 = 1$ and
\begin{equation*}
    (\alpha)_k =\frac{\Gamma(\alpha + k)}{\Gamma(\alpha)} =\alpha(\alpha +1)...(\alpha + k -1),
\end{equation*}
for $k = 1,2,...$
\end{theorem}
\begin{proof}
If $B_t = \frac{1}{2}Z_{4t}$, $S_t=e^{B_t-t}$, then $Y_{4t}= S^2_t$.
So
\begin{align*}
dS_t = S_t(dB_t - \frac{1}{2}dt), \qquad S_0 =1,
\end{align*}
 and
\begin{align} \label{eq-przed-Y}
  \Gamma_{4t} =  \frac{Y_{4t}}{1+\beta\int_0^{4t}Y_sds} =
    \frac{S^2_t}{1+4\beta\int_0^tS^2_u du}.
\end{align}
Define a new process
\begin{equation} \label{def:theta1}
    \theta_t = \frac{1}{4\beta} \frac{1+4\beta\int_0^t S^2_u du}{S^2_t}.
\end{equation}
From the It\^o lemma
\begin{align*}
d\theta_t = -2\theta_tdB_t + (4\theta_t +1)dt,  \qquad \theta_0 =
\frac{1}{4\beta}.
\end{align*}
Now we observe that the diffusion $\theta$ has the generator
\begin{equation}
    \mathcal{A}_{\theta} = 2x^2\frac{d^2}{dx^2} + (4x+1)\frac{d}{dx},
\end{equation}
which is identical with the generator of the process
$$\chi_t =
\exp(2B_t+2t)\Big(\frac{1}{4\beta} + \int_0^t\exp(-2B_u-2u)du\Big)$$ since
from the It\^o lemma
$$
    d\chi_t = 2\chi_tdB_t + (4\chi_t + 1)dt .
$$
Hence and from the fact that $\chi_0 = \theta_0$ we deduce that
processes $\theta$ and $\chi$ have the same distribution.
Let us take another Brownian motion $B_t^* =
\sqrt{2}B_{\frac{t}{2}}$ and define the process
\begin{align}
     \eta_t = \exp(\sqrt{2}B^*_t + t)\Big(\frac{1}{2\beta} +
     \int_0^t\exp(-\sqrt{2}B^*_u-u)du\Big).
\end{align}
Using the fact that $\theta_t$ and $\chi_t$ have the same
distribution, we obtain $ 2\theta_{\frac{t}{2}} =^d \eta_t$.
Moreover, we know that $\eta_t$ is a Markov process with the
resolvent
\begin{align}
    U_{\lambda}f(x) = \frac{\Gamma((1+\sqrt{4\lambda + 1})/2)}{\Gamma(1+\sqrt{4\lambda + 1})}\Big(\phi_1(x)\int_0^x e^{-\frac{1}{y}}\phi_2(y)f(y)dy+\phi_2(x)\int_x^{\infty}e^{-\frac{1}{y}}\phi_1(y)f(y)dy\Big).
\end{align}
(for details see   \cite[Theorem 3.1]{DGY}), so we conclude by
Theorem \ref{tw:momrep1}, \eqref{eq-przed-Y}, \eqref{def:theta1} and
definition of $\eta$ that
\begin{align*}
    \mathbb{E}X_{2T_{\lambda}}^{\alpha} =
     e^{-(\alpha\rho+\beta)}\mathbb{E}\exp\Big\{\frac{\alpha\rho+\beta}{2\beta}\frac{1}{\eta_{T_{\lambda}}}\Big\}
= \frac{1}{\lambda}
e^{-(\alpha\rho+\beta)}U_{\lambda}f\Big(\frac{1}{2\beta}\Big),
\end{align*}
with $f(x) =
\exp\Big\{\frac{\alpha\rho+\beta}{2\beta}\frac{1}{x}\Big\}$.
\end{proof}

\subsection{Moments  of the asset price in the  Stein model}

 In this section we consider the Stein model, i.e. the model
described by \eqref{defX}  with $Y$ being  an Ornstein-Uhlenbeck
process, so  $Y$ is  given by \eqref{defY-2}:
$$
  dY_t = -
\lambda Y_t dt + dZ_t,  \qquad Y_0 = 1, \quad \lambda > 0 . $$
 For
$t$ in neighborhood of zero we find an exact value of $
\mathbb{E}X_t^{\alpha}$. Let $b$ be the unique solution of equation
\begin{equation}  \label{def-gran}
b(1-e^{-2b})= 2 .
\end{equation}
 \begin{proposition} \label{prop:St}
  Let  $\alpha>0$ and $\rho$ be such that $\alpha (1 - \rho^2) < 1$,
$\rho < \lambda / \alpha$ and
 \begin{align} \label{df:gam}
    \gamma^2 &= \lambda^2 -\alpha^2(1-\rho^2)+\alpha - 2\lambda\rho\alpha > 0.
\end{align}
If $t \in [0, \frac{b}{\lambda}),$
then
\begin{equation}
    \mathbb{E}X_t^{\alpha} = e^{(1+t)\beta}\Big(\cosh(\gamma t)
    + \frac{2\beta}{\gamma}\sinh(\gamma t)\Big)^{-\frac{1}{2}},
\end{equation}
where $b$ is given by \eqref{def-gran} and
\begin{align*}
    \beta &= \frac{1}{2}(\lambda - \rho\alpha).
\end{align*}
\end{proposition}
\begin{proof}
Define the new measure
\begin{equation}
    \frac{d\mathbb{Q}}{d\mathbb{P}}\Big|_{\mathcal{F}_t} =
    e^{\lambda\int_0^tY_udZ_u -\frac{\lambda^2}{2}\int_0^tY_u^2du}.
\end{equation}
Clearly, $Y_t$ is the Gaussian random variable with the mean
$e^{-\lambda t}$ and variance $\frac{1}{2\lambda}(1-e^{-2\lambda
t})$.
 Moreover $\mathbb{E}e^{\frac{t\lambda^2}{2}Y_u^2}du <\infty$ for
$u<t$, since $\frac{\lambda t}{2}(1-e^{-2\lambda t}) < 1$ by
assumption on $t$. Thus, by Jensen inequality,
\begin{align*}
    \mathbb{E}e^{\frac{\lambda^2}{2}\int_0^tY_u^2du}&
    \leq \mathbb{E}
    \Big(\frac{1}{t}\int_0^te^{\frac{t\lambda^2}{2}Y_u^2}du \Big)
    =
    \frac{1}{t}\int_0^t\mathbb{E}e^{\frac{t\lambda^2}{2}Y_u^2}du \\
   &=
    \frac{1}{t}\int_0^t\frac{e^{e^{-2\lambda u}u\lambda^2 / (2-u\lambda(1-\exp (-2\lambda u)))}}{\sqrt{1-\frac{u\lambda}{2}(1-e^{-2\lambda u})}}\ du\\
    &< \frac{1}{t}\int_0^t\frac{e^{u\lambda^2 / 2}}{\sqrt{1-\frac{u\lambda}{2}(1-e^{-2\lambda u})}}\ du
    < \infty.
\end{align*}
So, $\mathbb{Q}$ is a probability measure, by the Novikov criterion.
Observe, by the Girsanov theorem,
 that  the process $Y$ is a  Brownian motion under $\mathbb{Q}$ starting from $1$.
Formula  \eqref{defY-2} implies
$$\int_0^tY_udZ_u = \frac12 (Y_t^2
-(t+1)) +\lambda\int_0^tY_u^2du, $$
 so by \eqref{potega} we obtain
\begin{align*}
\mathbb{E}X_t^{\alpha} &=
e^{-(t+1)\frac{\alpha\rho-\lambda}{2}}\mathbb{E}_{\mathbb{Q}}e^{\frac{\alpha\rho-\lambda}{2}Y_t^2+\frac{2\lambda\alpha
\rho -\lambda^2+\alpha^2(1-\rho^2)-\alpha}{2}\int_0^tY_u^2du}\\
&= e^{(1+t)\beta}\mathbb{E}e^{-\beta B_t^2 -\frac{\gamma^2}{2}\int_0^tB_u^2du}\\
&= e^{(1+t)\beta}\Big(\cosh(\gamma t) +
\frac{2\beta}{\gamma}\sinh(\gamma t)\Big)^{-\frac{1}{2}}
\exp\Big(\frac{\beta}{2} -\frac{(\gamma/2 + \beta)e^{\gamma t}}{\cosh(\gamma t) +
\frac{2\beta}{\gamma}\sinh(\gamma t)}\Big),
\end{align*}
 where in the last equation we use the form of the
Laplace'a transform of $(B_t^2,\int_0^tB_u^2du)$, where $B$ is a
 Brownian motion starting from $1$ (see e.g.  \cite[formula 1.9.7 page 168]{SB}).
\end{proof}
\begin{remark}
If $\alpha \in (0,1)$, then \eqref{df:gam} holds.
If $\alpha\geq 1$, $\alpha (1 - \rho^2) < 1 - 2\lambda\rho$ and $\rho\leq 0$ then \eqref{df:gam} holds.
\end{remark}

\bibliographystyle{plain}

\end{document}